# ANOTHER CORRECTION
## ERROR ESTIMATES FOR BINOMIAL APPROXIMATIONS OF GAME OPTIONS

*THE ANNALS OF APPLIED PROBABILITY* **16** (2006) 984–1033

BY YAN DOLINSKY AND YURI KIFER

*Hebrew University*

This note deals with the substantial inaccuracies in Lemmas 3.4 and 3.5 [more specifically, inequalities (4.48) and (4.53) of their proofs] and in Theorems 2.2 and 2.3 of [1] related to the important point that if a game option is not exercised or canceled before the expiration (horizon) time then the seller pays no penalty to the buyer, which is natural but does not agree well with the direct extension of payoff formulas beyond the horizon. The arguments in [1] do not require any modification if penalties in the corresponding game options are extended by zero beyond the horizon which, in view of the Lipschitz-type condition (2.2) there, would be a somewhat restrictive requirement since it eliminates the case of a constant (nonzero) penalty. Of course, there is no problem with the argument there if we consider just the American options case. We will deal with Theorems 2.2 and 2.3 later on in this note (warning the reader that our first correction of the proof there contains inaccuracies) and we start with showing that the estimate of Theorem 2.1 remains true if in place of Lemmas 3.4–3.6 we employ the following argument which extends the idea of Lemma 3.6 there. In the notations, the value of a game option in the Black–Scholes market is given by

$$(1) \qquad V(z) = \inf_{\sigma \in \mathcal{T}_{0T}^B} \sup_{\tau \in \mathcal{T}_{0T}^B} E^B Q_z^B(\sigma, \tau)$$

which in view of Lemma 3.3 in [1] should be compared with

$$(2) \qquad V_n^{B,\theta}(z) = \inf_{\zeta \in \mathcal{T}_{0,n}^{B,n}} \sup_{\eta \in \mathcal{T}_{0,n}^{B,n}} E^B Q_z^B(\theta_\zeta^{(n)}, \theta_\eta^{(n)}).$$

LEMMA 1. *There exists a constant $C > 0$ such that for all $z, n > 0$,*

(3) $$|V(z) - V_n^{B,\theta}(z)| \leq C(F_0(z) + \Delta_0(z) + z + 1)n^{-1/4}.$$

PROOF. For any $\delta > 0$ there exists $\sigma_\delta \in \mathcal{T}_{0T}^B$ such that

(4) $$V(z) \geq \sup_{\tau \in \mathcal{T}_{0T}^B} E^B Q_z^B(\sigma_\delta, \tau) - \delta.$$

As in [1] for each $\sigma \in \mathcal{T}_{0T}^B$ set $\nu_\sigma = \min\{k \in \mathbb{N} : \theta_k^{(n)} \geq \sigma\}$ and define

$$\alpha_{\sigma,\delta} = (n \wedge \nu_{\sigma_\delta})\mathbb{I}_{\sigma_\delta < T} + n\mathbb{I}_{\sigma_\delta = T}.$$

Since the Brownian filtration is right continuous it is easy to see that

$$\{\alpha_{\sigma,\delta} = n\} = \{\nu_{\sigma_\delta} \geq n\} \cup (\{\nu_{\sigma_\delta} < n\} \cap \{\sigma_\delta = T\})$$
$$= \{\theta_{n-1}^{(n)} < \sigma_\delta\} \cup (\{\theta_{n-1}^{(n)} \geq \sigma_\delta\} \cap \{\sigma_\delta = T\}) \in \mathcal{F}_{\theta_{n-1}^{(n)}}^B$$

and

$$\{\alpha_{\sigma,\delta} \leq k\} = \{\nu_{\sigma_\delta} \leq k\} \cap \{\sigma_\delta < T\} = \{\theta_k^{(n)} \geq \sigma_\delta\} \cap \{\sigma_\delta < T\} \in \mathcal{F}_{\theta_k^{(n)}}^B.$$

Hence $\alpha_{\sigma,\delta} \in \mathcal{T}_{0,n}^{B,n}$. Now, let $\beta_{\sigma,\delta} \in \mathcal{T}_{0,n}^{B,n}$ satisfies

(5) $$E^B Q_z^B(\theta_{\alpha_{\sigma,\delta}}^{(n)}, \theta_{\beta_{\sigma,\delta}}^{(n)}) \geq \sup_{\eta \in \mathcal{T}_{0,n}^{B,n}} E^B Q_z^B(\theta_{\alpha_{\sigma,\delta}}^{(n)}, \theta_\eta^{(n)}) - \delta.$$

Then by (2), (4) and (5),

(6) $$V_n^{B,\theta}(z) - V(z) \leq E^B Q_z^B(\theta_{\alpha_{\sigma,\delta}}^{(n)}, \theta_{\beta_{\sigma,\delta}}^{(n)}) - \sup_{\tau \in \mathcal{T}_{0T}^B} E^B Q_z^B(\sigma_\delta, \tau) + 2\delta$$
$$\leq E^B(Q_z^B(\theta_{\alpha_{\sigma,\delta}}^{(n)}, \theta_{\beta_{\sigma,\delta}}^{(n)}) - Q_z^B(\sigma_\delta, \theta_{\beta_{\sigma,\delta}}^{(n)} \wedge T)) + 2\delta.$$

It is clear that if $\theta_{\alpha_{\sigma,\delta}}^{(n)} < \theta_{\beta_{\sigma,\delta}}^{(n)}$, then $\sigma_\delta < \theta_{\beta_{\sigma,\delta}}^{(n)} \wedge T$,

(7) $$Q_z^B(\theta_{\alpha_{\sigma,\delta}}^{(n)}, \theta_{\beta_{\sigma,\delta}}^{(n)}) - Q_z^B(\sigma_\delta, \theta_{\beta_{\sigma,\delta}}^{(n)} \wedge T) \leq J_1 + J_2 + J_3,$$

where

$$J_1 = |F_{\theta_{\beta_{\sigma,\delta} \wedge \alpha_{\sigma,\delta}}^{(n)}}(S^B(z)) - F_{\theta_{\beta_{\sigma,\delta}}^{(n)} \wedge \sigma_\delta}(S^B(z))|,$$

$$J_2 = |G_{\theta_{\alpha_{\sigma,\delta}}^{(n)}}(S^B(z)) - G_{\sigma_\delta}(S^B(z))|$$

and

$$J_3 = |e^{-r\theta^{(n)}_{\beta_{\sigma,\delta}} \wedge \alpha_{\sigma,\delta}} - e^{-r(\theta^{(n)}_{\beta_{\sigma,\delta}} \wedge \sigma_\delta)}| F_{\theta^{(n)}_{\beta_{\sigma,\delta}}}(S^B(z))$$

$$+ |e^{-r\theta^{(n)}_{\alpha_{\sigma,\delta}}} - e^{-r\sigma_\delta}| G_{\sigma_\delta}(S^B(z)).$$

By assumption (2.2),

(8)
$$\begin{aligned} J_1 \leq L\bigg(\bigg(|\theta^{(n)}_n - T| + \max_{1\leq k\leq n}|\theta^{(n)}_k \wedge T - \theta^{(n)}_{k-1} \wedge T|\bigg) \\ \times \bigg(1 + ze^{r(\theta^{(n)}_n \vee T)} \sup_{0\leq t\leq \theta^{(n)}_n \vee T} e^{\kappa B_t}\bigg) \\ + \sup_{\theta^{(n)}_n \vee T \geq t \geq \theta^{(n)}_n \wedge T} |S^B_t(z) - S^B_{\theta^{(n)}_n \wedge T}(z)| \\ + \max_{1\leq k\leq n} \sup_{\theta^{(n)}_k \wedge T \geq t \geq \theta^{(n)}_{k-1} \wedge T} |S^B_t(z) - S^B_{\theta^{(n)}_{k-1}}(z)|\bigg), \end{aligned}$$

where $\theta^{(n)}_0 = 0$. Employing (4.7), (4.8) and (4.25), we obtain in the same way as in (4.60)–(4.67) that

(9) $$E^B J_1 \leq C^{(1)} z n^{-1/4}$$

for some $C^{(1)} > 0$ independent of $z$ and $n$. Next, by (2.2), the same bound holds for $J_2$, and so

(10) $$E^B J_2 \leq C^{(1)} z n^{-1/4}.$$

Finally, by (2.3) and (2.4),

(11)
$$\begin{aligned} J_3 \leq 2rL(\theta^{(n)}_n \vee T + 2)(F_0(z) + \Delta_0(z) + z)e^{r(\theta^{(n)}_n \vee T)} \\ \times \bigg(|\theta^{(n)}_n - T| + \max_{1\leq k\leq n}|\theta^{(n)}_k \wedge T - \theta^{(n)}_{k-1} \wedge T|\bigg)\bigg(1 + \sup_{0\leq t\leq \theta^{(n)}_n \vee T} e^{\kappa B_t}\bigg). \end{aligned}$$

Again, employing (4.7) and (4.8), we derive that

(12) $$E^B J_3 \leq C^{(2)}(F_0(z) + \Delta_0(z) + z) n^{-1/4}.$$

Now, (6), (7), (9), (10) and (12) yield that

(13) $$V^{B,\theta}_n(z) - V(z) \leq (2C^{(1)} z + C^{(2)}(F_0(z) + \Delta_0(z) + z)) n^{-1/4} + 2\delta.$$

In order to estimate this difference in the other direction choose $\zeta_\delta \in \mathcal{T}^{B,n}_{0,n}$ such that

(14) $$V^{B,\theta}_n(z) \geq \sup_{\eta \in \mathcal{T}^{B,n}_{0,n}} E^B Q^B_z(\theta^{(n)}_{\zeta_\delta}, \theta^{(n)}_\eta) - \delta.$$

Set
$$\sigma_\delta = T\mathbb{I}_{\zeta_\delta = n} + \theta^{(n)}_{\zeta_\delta} \wedge T\mathbb{I}_{\zeta_\delta < n}.$$

For any $t \leq T$,
$$\{\sigma_\delta < t\} = \{\zeta_\delta < n\} \cap \{\theta^{(n)}_{\zeta_\delta} < t\} \in \mathcal{F}^B_t$$

and
$$\{\sigma_\delta = T\} = \Omega \setminus \{\sigma_\delta < T\} \in \mathcal{F}^B_T.$$

Since the Brownian filtration is right continuous we conclude that $\sigma_\delta \in \mathcal{T}^B_{0T}$. Let $\tau_\delta \in \mathcal{T}^B_{0T}$ satisfies
$$E^B Q^B_z(\sigma_\delta, \tau_\delta) + \delta \geq \sup_{\tau \in \mathcal{T}^B_{0T}} E^B Q^B_z(\sigma_\delta, \tau)$$

then by (1) and (14) we obtain

(15) $\quad V(z) - V^{B,\theta}_n(z) \leq 2\delta + E^B(Q^B_z(\sigma_\delta, \tau_\delta) - Q^B_z(\theta^{(n)}_{\zeta_\delta}, \theta^{(n)}_{\nu_{\tau_\delta} \wedge n})),$

where, recall, $\nu_\tau = \min\{k \in \mathbb{N} : \theta^{(n)}_k \geq \tau\}$. As before, we obtain that $\nu_{\tau_\delta} \wedge n \in \mathcal{T}^{B,n}_{0,n}$.

It is easy to see from the definition that $\sigma_\delta < \tau_\delta$ if and only if $\theta^{(n)}_{\zeta_\delta} < \theta^{(n)}_{\nu_{\tau_\delta} \wedge n}$, and so

(16) $\quad Q^B_z(\sigma_\delta, \tau_\delta) - Q^B_z(\theta^{(n)}_{\zeta_\delta}, \theta^{(n)}_{\nu_{\tau_\delta} \wedge n}) \leq J_4 + J_5 + J_6,$

where
$$J_4 = |F_{\tau_\delta}(S^B(z)) - F_{\theta^{(n)}_{\nu_{\tau_\delta} \wedge n}}(S^B(z))|,$$
$$J_5 = |G_{\sigma_\delta}(S^B(z)) - G_{\theta^{(n)}_{\zeta_\delta}}(S^B(z))|$$

and
$$J_6 = |e^{-r\tau_\delta} - e^{-r\theta^{(n)}_{\nu_{\tau_\delta} \wedge n}}|F_{\tau_\delta}(S^B(z))$$
$$+ |e^{-r\sigma_\delta} - e^{-r\theta^{(n)}_{\zeta_\delta}}|G_{\sigma_\delta}(S^B(z)).$$

By (2.2),
$$J_4 \leq L\bigg(\bigg(|\theta^{(n)}_n - T| + \max_{1 \leq k \leq n}|\theta^{(n)}_k \wedge T - \theta^{(n)}_{k-1} \wedge T|\bigg)$$
$$\times \bigg(1 + ze^{r(\theta^{(n)}_n \vee T)} \sup_{0 \leq t \leq \theta^{(n)}_n \vee T} e^{\kappa B_t}\bigg)$$
(17)

$$+ \sup_{\theta_n^{(n)} \vee T \geq t \geq \theta_n^{(n)} \wedge T} |S_t^B(z) - S_{\theta_n^{(n)} \wedge T}^B(z)|$$

$$+ \max_{0 \leq k \leq n} \sup_{\theta_k^{(n)} \wedge T \geq t \geq \theta_{k-1}^{(n)} \wedge T} |S_t^B(z) - S_{\theta_{k-1}^{(n)} \wedge T}^B(z)|\bigg).$$

Thus, $J_4$ has the same bound as $J_1$, and so

(18) $$E^B J_4 \leq C^{(1)} z n^{-1/4}.$$

Next, by (2.2) and the definition of $\sigma_\delta$,

(19)
$$J_5 \leq L\bigg(|\theta_n^{(n)} - T|\bigg(1 + ze^{r(\theta_n^{(n)} \vee T)} \sup_{0 \leq t \leq \theta_n^{(n)} \vee T} e^{\kappa B_t}\bigg)$$
$$+ \sup_{\theta_n^{(n)} \vee T \geq t \geq \theta_n^{(n)} \wedge T} |S_t^B(z) - S_{\theta_n^{(n)}(z) \wedge T}^B(z)|$$
$$+ \max_{1 \leq k \leq n} \sup_{\theta_k^{(n)} \wedge T \geq t \geq \theta_{k-1}^{(n)} \wedge T} |S_t^B(z) - S_{\theta_{k-1}^{(n)}}^B(z)|\bigg).$$

Hence, $J_5$ can be estimated by the right-hand side of (8), and so

(20) $$E^B J_5 \leq C^{(1)} z n^{-1/4}.$$

Now, by (2.3) and (2.4),

(21)
$$J_6 \leq 2rL(\theta_n^{(n)} \vee T + 2)(F_0(z) + \Delta_0(z) + z)e^{rT}$$
$$\times \bigg(|\theta_n^{(n)} - T| + \max_{1 \leq k \leq n} |\theta_k^{(n)} \wedge T - \theta_{k-1}^{(n)} \wedge T|\bigg)$$
$$\times \bigg(1 + \sup_{0 \leq t \leq \theta_n^{(n)} \vee T} e^{\kappa B_t}\bigg)$$

and in the same way as in (12) we obtain that

(22) $$E^B J_6 \leq C^{(2)}(F_0(z) + \Delta_0(z) + z)n^{-1/4}.$$

Since $\delta > 0$ is arbitrary we obtain (3) from (13), (15), (16), (18), (20) and (22).

Next we deal with Theorems 2.2 and 2.3. First, we have to replace the stopping time $\varphi_n^* = \theta_{\mu_n^* \circ \lambda_B^{(n)}}^{(n)}$ in Theorem 2.2 by $\sigma^* = (T \wedge \varphi_n^*)\mathbb{I}_{\mu_n^* \circ \lambda_B^{(n)} < n} + T\mathbb{I}_{\mu_n^* \circ \lambda_B^{(n)} = n}$ and the stopping time $\theta_\varphi^{(n)}$ in Theorem 2.3 by $\sigma = (T \wedge \theta_\varphi^{(n)})\mathbb{I}_{\varphi < n} + T\mathbb{I}_{\varphi = n}$ where, recall, $\varphi = \mu \circ \lambda_B^{(n)}$ and $(\mu \circ \lambda_\xi^{(n)}, \pi)$ is a hedge for some self-financing portfolio strategy $\pi$ in the corresponding binomial market. Then the proof of Theorem 2.2 should be modified according to the arguments below which will give (2.17) with $\sigma^*$ in place of $\varphi_n^*$ and the same $\psi_n^*$ as there. In place of (2.22) in Theorem 2.3 we will have the following shortfall estimates:

$$\sup_{\tau \in \mathcal{T}_{0T}^B} E^B(R_z^B(\sigma,\tau) - Z_{\sigma \wedge \tau}^B)^+ \leq C(F_0(z) + \Delta_0(z) + z + 1)n^{-1/4}(\ln n)^{3/4} \quad (23)$$

for some constant $C > 0$ and, furthermore, for any $\epsilon > 0$ there exists a constant $C_\epsilon$ such that

$$(24) \quad E^B \sup_{0 \leq t \leq T} E^B(R_z^B(\sigma,t) - Z_{\sigma \wedge t}^B)^+ \leq C_\epsilon(F_0(z) + \Delta_0(z) + z + 1)n^{-1/4-\epsilon}.$$

Here $Z_t^B$ is the portfolio value at time $t$ for a self-financing strategy $\pi^B$ in the Black–Scholes market constructed by $\pi$ in Theorem 2.3 and we have to specify that between the times $T \wedge \theta_n^{(n)}$ and $T$ we manage $Z_t^B$ so that its discounted value $\check{Z}_t^B$ stays constant, that is, $\check{Z}_t^B = \check{Z}_{T \wedge \theta_n^{(n)}}^B$ for all $t \in [T \wedge \theta_n^{(n)}, T]$. This is done by selling all stocks at the time $\theta_n^{(n)}$ if $\theta_n^{(n)} < T$ (otherwise, doing nothing), buying immediately bonds for all money and doing nothing afterward.

In order to prove (23) and (24) set

$$\Psi = \sup_{0 \leq t \leq T} \left| Q_z^B(\sigma, t) - Q_z^{B,n}\left(\frac{\varphi T}{n}, \frac{(n \wedge \nu_t)T}{n}\right) \right|,$$

where, recall, $Q_z^B(s,t) = e^{-rs \wedge t} R_z^B(s,t)$ and $Q_z^{B,n}(s,t) = e^{-rs \wedge t} R_z^{B,n}(s,t)$ are the discounted payoffs. By (4.22)–(4.28) it follows that

$$(25) \quad \max_{0 \leq k \leq n} \max_{0 \leq l \leq n} \left| Q_z^{B,n}\left(\frac{kT}{n}, \frac{lT}{n}\right) - Q_z^{B,\theta,n}(\theta_k^{(n)}, \theta_l^{(n)}) \right|$$
$$\leq r \max_{0 \leq k \leq n} \left| \theta_k^{(n)} - \frac{kT}{n} \right| \max_{0 \leq k \leq n} G_{kT/n}(S^{B,n}(z)) + 3I,$$

where $I$ is the right-hand side of (4.28). As in Lemma 3.2 we obtain

$$(26) \quad E^B \max_{0 \leq k \leq n} \max_{0 \leq l \leq n} E^B \left| Q_z^{B,n}\left(\frac{kT}{n}, \frac{lT}{n}\right) - Q_z^{B,\theta,n}(\theta_k^{(n)}, \theta_l^{(n)}) \right|$$
$$\leq C^{(1)}(F_0(z) + \Delta_0(z) + z + 1)n^{-1/4}(\ln n)^{3/4}.$$

Similarly to (16)–(22) we conclude that

$$(27) \quad E^B \sup_{0 \leq t \leq T} |Q_z^B(\sigma,t) - Q_z^B(\theta_\varphi^{(n)}, \theta_{n \wedge \nu_t}^{(n)})| \leq C^{(2)}(F_0(z) + \Delta_0(z) + z)n^{-1/4}$$

which together with (26) and Lemma 3.3 yields that

$$(28) \quad E^B \Psi \leq C^{(3)}(F_0(z) + \Delta_0(z) + z + 1)n^{-1/4}(\ln n)^{3/4}.$$

From (5.24) it follows that $\check{Z}^B_{\varphi\wedge\nu_\tau} \geq Q^{B,n}_z(\frac{T\varphi}{n}, \frac{T(\nu_\tau\wedge n)}{n})$ for any $\tau \in \mathcal{T}^B_{0T}$ where, recall, that $\check{Z}^B_t = e^{-rt}Z^B_t$ is the discounted portfolio value. Since $\check{Z}^B_t$ is a martingale which does not change on the time interval $[T\wedge\theta^{(n)}_n, T]$ and since $\sigma\wedge\tau\wedge\theta^{(n)}_n \leq \theta^{(n)}_{\nu_\tau\wedge\varphi}$ we obtain

$$\check{Z}^B_{\sigma\wedge\tau} = \check{Z}^B_{\sigma\wedge\tau\wedge\theta^{(n)}_n} = E^B(\check{Z}^B_{\theta^{(n)}_{\nu_\tau\wedge\varphi}}|\mathcal{F}^B_{\sigma\wedge\tau\wedge\theta^{(n)}_n}) = E^B(\check{Z}^B_{\theta^{(n)}_{\nu_\tau\wedge\varphi}}|\mathcal{F}^B_{\sigma\wedge\tau})$$
(29)
$$\geq E^B(Q^B_z(\sigma,\tau) - \Psi|\mathcal{F}^B_{\sigma\wedge\tau}) = Q^B_z(\sigma,\tau) - E^B(\Psi|\mathcal{F}^B_{\sigma\wedge\tau})$$

taking into account that $\check{Z}^B_{\theta^{(n)}_{\nu_\tau\wedge\varphi}}$ is $\mathcal{F}^B_{\theta^{(n)}_n}$-measurable and $Q^B_z(\sigma,\tau)$ is $\mathcal{F}^B_{\sigma\wedge\tau}$-measurable. We use the latter measurability in the last equality and the former one in the third equality of (29) together with the formula $E^B(X|\mathcal{F}^B_{\tilde\sigma\wedge\tilde\tau}) = E^B(X|\mathcal{F}^B_{\tilde\tau})$ which is valid provided $X$ is $\mathcal{F}^B_{\tilde\tau}$-measurable and $\tilde\tau, \tilde\sigma$ are stopping times. This together with (28) gives

$$\sup_{\tau\in\mathcal{T}^B_{0,T}} E^B(R^B_z(\sigma,\tau) - Z^B_{\sigma\wedge\tau})^+ \leq e^{rT}\sup_{\tau\in\mathcal{T}^B_{0,T}} E^B(Q^B_z(\sigma,\tau) - \check{Z}^B_{\sigma\wedge\tau})^+$$

$$\leq e^{rT}\sup_{\tau\in\mathcal{T}^B_{0,T}} E^B(E^B(\Psi|\mathcal{F}^B_{\sigma\wedge\tau})) = e^{rT}E^B(\Psi)$$

$$\leq C^{(4)}(F_0(z) + \Delta_0(z) + z + 1)n^{-1/4}(\ln n)^{3/4}$$

proving (23).

Next, fix $\alpha > 0$ and set $\delta = \frac{\sqrt{1+\alpha}-1}{2}$. By (2.3) and (4.8), together with the Hölder inequality and (28), we obtain (for sufficiently large $n$) that

(30)
$$E^B\Psi^{1+\delta} \leq E^B\left(\Psi^{1/(1+\delta)}\left(\sup_{0\leq t\leq T} G_t(S^B(z))\right)^{1+\delta-1/(1+\delta)}\right)$$
$$\leq C_\delta(E^B\Psi)^{1/(1+\delta)} \leq \tilde{C}_\delta(F_0(z) + \Delta_0(z) + z + 1)n^{-1/(4(1+2\delta))}.$$

Observe that $\{E^B(\Psi|\mathcal{F}^B_{\sigma\wedge t})\}^T_{t=0}$ is a regular martingale and so from (29), (30) and the Doob's maximal inequality we obtain

$$E^B \sup_{0\leq t\leq T}((R^B_z(\sigma,t) - Z^B_{\sigma\wedge t})^+)^{1+\delta}$$

$$\leq e^{rT(1+\delta)}E^B \sup_{0\leq t\leq T}((Q^B_z(\sigma,t) - \check{Z}^B_{\sigma\wedge t})^+)^{1+\delta}$$

$$\leq e^{rT(1+\delta)}E^B \sup_{0\leq t\leq T}(E^B(\Psi|\mathcal{F}^B_{\sigma\wedge t}))^{1+\delta} \leq C^{(1)}_\delta E^B\Psi^{1+\delta}$$

$$\leq C^{(2)}_\delta(F_0(z) + \Delta_0(z) + z + 1)n^{-1/(4(1+2\delta))},$$

and using Jensen's inequality we estimate the left-hand side of (24) by $C_\alpha n^{-1/(4(1+\alpha))}$, completing the proof. $\square$

INSTITUTE OF MATHEMATICS
HEBREW UNIVERSITY
JERUSALEM 91904
ISRAEL
E-MAIL: yann1@math.huji.ac.il
         kifer@math.huji.ac.il